\documentclass[a4paper]{amsproc}

\newtheorem{theorem}{Theorem}[section]

\theoremstyle{definition}
\newtheorem{definition}[theorem]{Definition}

\theoremstyle{remark}

\numberwithin{equation}{section}

\usepackage[all]{xy}

\DeclareMathOperator{\ind}{ind}
\DeclareMathOperator{\interior}{int}
\DeclareMathOperator{\im}{im}
\DeclareMathOperator{\id}{id}

\begin{document}

\title{Bordism Invariance of the Coarse Index}


\author{Christopher Wulff}

\subjclass[2000]{Primary 19K56, 58J22; Secondary 19K35}

\keywords{}

\date{}

\dedicatory{}

\begin{abstract}
We prove bordism invariance of the coarse index of complex elliptic
pseudodifferential operators.
In our discussion we introduce directed $c$-bordisms, whose usefulness is illustrated in the context of existence of uniformly positive scalar curvature metrics on open manifolds.
\end{abstract}

\maketitle

\section{Introduction}
One way to generalize the analytic index of elliptic pseudodifferential operators to
non-compact manifolds is the coarse index as defined in \cite[Chapter 12]{MR1817560}. It is defined for any manifold equipped with a
proper coarse structure (see \cite[Chapter 6]{MR1817560} or
\cite{MR2007488} as
references on coarse geometry).
The coarse index of the spinor Dirac operator is an obstruction to the existence of
uniformly positive scalar curvature much as the analytic index in the compact case \cite[Chapter 12]{MR1817560}.

In index theory, bordism invariance has always been an important
tool. Ever since Atiyah and Singer established it in the first proof of their index theorem \cite[Chapter XVII]{MR0198494}, a variety of new proofs for different settings in compact and other Fredholm operator cases have been published, for example
\cite{MR1113688,MR1231967,MR1402879,MR1873784,MR2139923,MR2274155,MR2231925}. In the present paper we discuss bordism invariance of the coarse index.

An essential aspect is to find an appropriate notion of bordism in the coarse geometric setting. A discussion which we will postpone until section \ref{cbordisms} suggests the following definition:
\begin{definition}\label{cbordism}
A \emph{directed $c$-bordism} is a triple $(W;N_1,N_2)$ of smooth
manifolds  such that
\begin{enumerate}
\item the boundary of $W$ decomposes as $\partial W=N_1\dot\cup N_2$, \item each of $W,N_1,N_2$ is equipped with a proper coarse structure, such that
the inclusion $\iota_1:N_1\xrightarrow{\subset}W$ is a coarse map and the inclusion
$\iota_2:N_2\xrightarrow{\subset}W$ is a coarse equivalence.
\end{enumerate}
\end{definition}

Our notion of operator bordism will be $K$-homological:
Graded (resp. ungraded) symmetric elliptic pseudodifferential operators $D_M$ over a manifold $M$ determine $K$-homology classes $[D_M]\in K_0(M)$ (resp. $\in K_1(M)$) \cite[Chapter 10]{MR1817560} (see also
\cite{MR1153932} for the pseudodifferential case). We consider graded ($p=0$) or ungraded ($p=1$) operators $D_1,D_2$ over $N_1,N_2$ as bordant by a directed $c$-bordism $(W;N_1,N_2)$ if there exists a class $x\in K_{1-p}(\interior W)$ ($\interior W=W\setminus N_1\cup N_2$) which is mapped to
$[D_1]$, resp. $-[D_2]$ by the connecting homomorphisms
$\partial_i:K_{1-p}(\interior W)\to K_p(N_i)$ of the six term exact sequences in $K$-homology for the pairs $(\interior W\cup N_i,N_i)$.

The definitions of operator bordisms in the literature on bordism invariance of the index can be subdivided into two (essentially equivalent) types: Topological ways of expressing operator bordisms in terms of the principal symbols of the operators or their symbol classes in $K$-theory are found in
\cite{MR0198494,MR2274155,MR2231925}.
The notion of symbol cobordism of \cite{MR2274155}, for example,
corresponds to our definition in the case $p=0$ under a natural
isomorphism
\begin{equation*}
K_*(M)\cong K^*(T^*M)
\end{equation*}
given in \cite{MR0266247,MR1153932}, which maps the $K$-homology class of an operator to its symbol class. This follows in the compact case from \cite[Proposition 6.9]{MR1153932}.

The other type of operators considered are graded or ungraded (mostly Dirac or first order differential) symmetric operators over the boundaries constructed analytically from ungraded resp. graded operators over the bordism as in \cite[Chapter 4]{MR1399087}. This setting is seen in \cite{MR1113688,MR2139923,MR1231967,MR1873784,MR1402879} in the case $p=0$.
Under minor simplifying assumptions, a direct calculation with the methods of \cite[Chapters 9,10]{MR1817560} shows that this analytic construction also gives rise to operators bordant in our sense. 

Our main theorem -- a directed version of bordism invariance -- now reads as follows:
\begin{theorem}\label{maintheorem}
For operators $D_1,D_2$ bordant by a directed $c$-bordism $(W;N_1,N_2)$, the coarse index of $D_1$ is mapped to the coarse index of $D_2$ by the canonical map 
\begin{equation}\label{Randabbildung}
K_p(C^*(N_1))\to K_p(C^*(W))\cong K_p(C^*(N_2))
\end{equation}
\end{theorem}
In the symmetric case where the inclusion of $N_1$ is a coarse equivalence as well, the map \eqref{Randabbildung} becomes a natural identification $K_p(C^*(N_1))\cong K_p(C^*(N_2))$.

The main ingredient to our proof of theorem \ref{maintheorem} is the naturality of the coarse assembly map $A_Y:K_*(Y)\to K_*(C^*(Y))$, which is defined for any locally compact space equiped with a proper coarse structure, under continuous coarse maps (see \cite[Chapter 12]{MR1817560} and \cite{MR1451755}). The coarse map $N_1\to N_2$ inducing
\ref{Randabbildung} however lacks continuity and hence does not induce a map of topological $K$-homology groups. This deficiency is resolved by using a bordism relating the two boundaries.

As an application we discuss in section \ref{application} quasi-isometry invariance of the coarse index in a special case, following an idea of Roe: Two copies of a connected manifold $M$ equiped with two different quasi-isometric complete Riemannian metrics and coarse structures
determined by the corresponding path metrics can be connected by a product $c$-bordism. The same
construction applies when the path metrics $d_0,d_1$ associated to the Riemannian metrics satisfy only $d_0\geq Cd_1-C'$ for some $C,C'>0$, yielding a directed product $c$-bordism together with a canonical map \begin{equation*}
K_p(C^*(M,d_0))\to K_p(C^*(M,d_1)).
\end{equation*}
This will be applied to the spinor Dirac operator in the context of existence of uniformly positive scalar curvature in Theorem
\ref{strengthening}.

\section{$c$-Bordisms}\label{cbordisms}
In this section we discuss the necessity and reasonability of various aspects of Definition \ref{cbordism}.
The definition of the coarse
indices $\ind_c(D_{1,2})$ of $D_{1,2}$
requires $N_{1,2}$ to be equipped with proper coarse structures.  This is usually taken to
be the metric coarse structure of the path metric on a complete connected Riemannian manifold.
Note however, that the coarse index is defined for any proper coarse structure, and indeed no relation between coarse structure and Riemannian metric, if present, is needed for our bordism invariance.

Note that any manifold $M$ is the boundary of a non-compact manifold $W$. In the special cases $W=M\times [0,t)$ equiped with a product coarse structure ($t\in(0,\infty]$ and $[0,t)$ equiped with the canonical metric coarse structure) the connecting homomorphism
$\partial_M:K_{1-p}(W\setminus M)\to K_p(M)$ of the pair $(W,M)$ is an isomorphism and thus its image contains $[D_M]$ for any operator $D_M$ over $M$.
Therefore, to prevent arbitrary operators from being null-bordant, the following restrictions to $W$ seem reasonable:
\begin{enumerate}
\item $W$ should carry a proper coarse structure as well, excluding examples as those above with $t<\infty$.
\item The inclusion of a part of $\partial W$ into $W$ should be a coarse equivalence, excluding cases like the one above with $t=\infty$. In definition \ref{cbordism} this part is $N_2$. \label{semiinfiniteproduct} \end{enumerate}

Possible weakenings of requirement \ref{semiinfiniteproduct} are covered by the following theorem describing the vanishing of the index of such null-bordant operators:

\begin{theorem}\label{vanishing}
Let $W$ be a manifold with boundary $M$ such that both $W,M$ carry proper coarse structures and $\iota:M\xrightarrow{\subset} W$ is a coarse map. If the $K$-homology class $[D_M]$ of an operator over $M$ lies in the image of the connecting homomorphism $\partial_M:K_{1-p}(W\setminus M)\to K_p(M)$, then
$$\ind_c(D_M)\in\ker(\iota_*:K_p(C^*M)\to K_p(C^*W)).$$
\end{theorem}
We will prove this theorem in the next section.

\section{Proof of bordism invariance}\label{invariance}
The coarse index $\ind_c(D_M)$ of a graded or ungraded symmetric elliptic pseudodifferential operator $D_M$ over a manifold $M$ equipped with a proper coarse structure is defined by applying the coarse assembly map $A_M$ to its $K$-homology class \cite[Chapter 12]{MR1817560}. Theorem \ref{vanishing} is a simple consequence of the naturality of the assembly map:
\begin{equation*}\xymatrix{
K_{1-p}(W\setminus M)\ar[r]^{\partial_M}
&K_p(M)\ar[r]^{\iota_*}\ar[d]_{A_M}
&K_p(W)\ar[d]^{A_W}
\\&K_p(C^*(M))\ar[r]_{\iota_*}
&K_p(C^*(W))
}\end{equation*}
commutes and the claim follows from exactness of the first row and $[D_M]\in\im(\partial_M)$.

Theorem \ref{maintheorem} can now be deduced from this proof. Given a directed $c$-bordism $(W;N_1,N_2)$, we equip $M=\partial W$ with the restricted
coarse structure and obtain a diagram of inclusions
\begin{equation*}\xymatrix{N_1\ar[r]^{i_1}\ar[dr]_{\iota_1}
&M\ar[d]^{\iota}
&N_2\ar[l]_{i_2}\ar[dl]^{\iota_2}
\\&W
}\end{equation*}
which are all coarse and continuous maps.
For the element $x\in K_{1-p}(\interior W)$ with $[D_1]=\partial_1x, [D_2]=-\partial_2x$ we have $\iota_*A_M\partial_Mx=0$ as above. By naturality of the exact sequence in $K$-homology for $C^*$-algebras under the inclusions 
\begin{equation*}
\kappa_i:(C_0(\interior W\cup N_i),C_0(N_i)) \xrightarrow{\subset} (C_0(W),C_0(M))
\end{equation*}
we see
\begin{equation*}
\kappa_i^*\partial_Mx=\partial_ix=(-1)^{i-1} [D_i]
\end{equation*}
and
\begin{equation*}
\kappa_1i_1^*+\kappa_2i_2^*=\id
\end{equation*}
on $C_0(M)$ implies
\begin{equation*}
(i_1)_*\kappa_1^*+(i_2)_*\kappa_2^*=\id
\end{equation*}
on $K_*(M)$. Therefore,
\begin{equation*}
\iota_*A_M(i_1)_*[D_1]+\iota_*A_M(i_2)_*(-[D_2])=0.
\end{equation*}
Using the naturality of the assembly map again we obtain
\begin{equation*}
(\iota_1)_*A_{N_1}[D_1]=(\iota_2)_*A_{N_2}[D_2].
\end{equation*}
This proves Theorem \ref{maintheorem}.

\section{Quasi-isometry invariance and directed
bordisms}\label{application}
A special case of quasi-isometry may be arranged as product $c$-bordism: Consider a smooth manifold $M$ carrying two distinct complete
quasi-isometric Riemannian metrics $g_0,g_1$ and coarse structures determined by the corresponding path metrics $d_0,d_1$. Choose a smooth function $\chi:[0,1]\to [0,1]$ with $\chi\equiv 0$ near $0$, $\chi\equiv 1$ near $1$ and equip $W=M\times [0,1]$ with Riemannian metric
$g=(1-\chi(t))g_0+\chi(t)g_1+dt^2$ (for the analytic construction of bordant dirac operators, a product metric on a collar neighbourhood of the boundary is useful) and coarse structure determined by the corresponding path metric. The inclusions of $M_i\equiv M\times
\{i\}=(M,g_i), i\in\{0,1\}$ as boundaries into $(W,g)$ are coarse
equivalences and we obtain a $c$-bordism $(W;M_0,M_1)$.

The same construction applies when the metrics are not quasi-isometric but satisfy only $d_0\geq Cd_1-C'$ for some $C,C'>0$, yielding a
\emph{directed} $c$-bordism with canonical map $K_*(M_0)\to K_*(M_1)$.

We illustrate the usefulness of this idea in the context of existence of uniformly positive scalar curvature on spin-manifolds: A spin structure on $M$ determines a spin structure on $W$ and the spinor Dirac operators $D_0,D_1$ of $M_0,M_1$ are bordant by the
fundamental class $x=[W,g]\in K_{\dim W}(\interior W)$ of $(W,g)$. Thus by bordism invariance, if the coarse index of $D_1$ does not vanish, neither does the coarse index of $D_0$ and therefore $g_0$ is not of uniformly positive scalar curvature.

As an application we revisit a result in \cite{MR2482205} concerning uniformly positive scalar curvature metrics on open manifolds. Consider a uniformly enlargeable manifold $M$ whose universal covering $\tilde M$ is spin, cf. \cite[Section 1]{MR2482205}.
In the proof of \cite[Corollary 1.8]{MR2482205} it is shown that the coarse index of the spinor Dirac operator of $(\tilde M,\tilde g_1)$ is nonzero for any metric $\tilde g_1$ on $\tilde M$ pulled back from a metric on $M$. Quasi-isometry invariance of the coarse index is employed to conclude that $\tilde M$ admits no Riemannian metric $\tilde g_0$ of uniformly positive scalar curvature quasi-isometric to $\tilde g_1$. Our directed product $c$-bordism allows the following strengthening:

\begin{theorem}\label{strengthening}
If $M$ is universally enlargeable and its universal covering $\tilde M$ is spin,
then $\tilde M$ does not admit a metric $\tilde g_0$ of uniformly positive scalar
curvature satisfying $\tilde d_0\geq C\tilde d_1-C'$ with $C,C'>0$ and $\tilde d_{0,1}$ being the path metrics of $\tilde g_{0,1}$ where $\tilde g_1$ is a metric pulled back from $M$.
\end{theorem}

\section*{Acknowledgements}
The present paper emanates from the author's diploma thesis at the Ludwig--Maximilians--Universit\"{a}t M\"{u}nchen, 2010, under the guidance of Bernhard Hanke.

\providecommand{\bysame}{\leavevmode\hbox to3em{\hrulefill}\thinspace}
\providecommand{\MR}{\relax\ifhmode\unskip\space\fi MR }
\providecommand{\MRhref}[2]{%
  \href{http://www.ams.org/mathscinet-getitem?mr=#1}{#2}
}
\providecommand{\href}[2]{#2}

\end{document}